\documentclass[10pt]{amsart}
\usepackage{url,graphicx}
\usepackage[dvips]{epsfig}
\usepackage{latexsym,color}
\usepackage{amscd,amssymb,verbatim}
\newcommand{\be}{\begin{enumerate}}
\newcommand{\ee}{\end{enumerate}}

\newcommand{\cc}{{\mathcal C}}
\newcommand{\cd}{{\mathcal D}}
\newcommand{\cf}{{\mathcal F}}

\newcommand{\car}{{\mathcal R}}

\newcommand{\bo}{\partial}
\newcommand{\bu}{\bullet}

\newcommand{\chom}{\text{\tt Hom}_{\,0}}
\newcommand{\cmp}{\complement}
\newcommand{\cs}{{\mathcal S}}

\newcommand{\dz}{{\mathbb Z}}

\newcommand{\im}{\text{\rm Im}\,}
\newcommand{\ind}{\text{\rm Ind}\,}

\newcommand{\join}[1]{\ast_{#1}}

\newcommand{\lra}{\longrightarrow}
\newcommand{\nin}{\noindent}
\newcommand{\pr}{\noindent{\bf Proof. }}

\newcommand{\ra}{\rightarrow}
\newcommand{\sgn}{\text{sgn}\,}
\newcommand{\sm}{\setminus}

\newcommand{\supp}{\text{\rm supp}\,}
\newcommand{\thin}{\text{\rm Thin}}
\newcommand{\thom}{\text{\tt Hom}\,}
\newcommand{\thomp}{\text{\tt Hom}_+}
\newcommand{\ti}{\tilde}

\newcommand{\wti}{\widetilde}
\newcommand{\xor}{\text{\tt xor}}
\newcommand{\zz}{{\dz_2}}

\newtheorem{thm}{Theorem}[section]
\newtheorem{df}[thm]{Definition}

\newtheorem{crl}[thm]{Corollary}
\newtheorem{prop}[thm]{Proposition}

\numberwithin{equation}{section}
\numberwithin{figure}{section}
\numberwithin{table}{section}

\begin{document}

\title
{Cohomology of colorings of cycles}

\author{Dmitry N. Kozlov}
\address{Department of Computer Science, Eidgen\"ossische Technische
Hochschule, Z\"urich, Switzerland}

\email{dkozlov@inf.ethz.ch}
\thanks {Research supported by Swiss National Science Foundation Grant PP002-102738/1}

\date\today

\begin{abstract}
We compute the cohomology groups of the spaces of colorings of cycles,
i.e., of the prodsimplicial complexes $\thom(C_m,K_n)$. We perform the
computation first with $\zz$, and then with integer coefficients. The
main technical tool is to use spectral sequences in conjunction with
a~detailed combinatorial analysis of a~family of cubical complexes,
which we call {\it torus front complexes}.

As an~application of our method, we demonstrate how to collapse each
connected component of $\thom(C_m,C_n)$ onto a~garland of cubes.
\end{abstract}

\maketitle

\section{Introduction}

Given a~graph $G$, one can try to color it with $n$ colors so that if
two vertices are connected by an~edge, then they should get different
colors. Traditionally a~question of interest has been ``does such
a~coloring exist?'' If it does, then the possible colorings will form
a~discrete set $\chom(G,K_n)$, with no a~priori given structure. 

It has recently been realized, that it is possible to add
an~additional topological structure to this set. This is done by
considering a~polyhedral complex $\thom(G,K_n)$, whose vertices are
all legal colorings, as the choice of our notations indicate, and
whose higher-dimensional cells encode ``commuting relations'' among
the colorings. The complex $\thom(G,K_n)$ is a~special case of the
{\it Lov\'asz Hom construction}, which associates to any pair of graphs
$T$ and $G$ (no double edges, loops are allowed) a~polyhedral, in fact
prodsimplicial, see \cite{IAS}, complex $\thom(T,G)$, whose vertices
are all graph homomorphisms from $T$ to $G$ and whose cells encode
``commuting relations'' among these homomorphisms.

This paper is devoted to studying the cohomology groups of the
coloring space of a~cycle: $\thom(C_m,K_n)$. When $m$ is an odd
number, the partial analysis of the topology of these spaces played
a~central role in the resolution of the Lov\'asz Conjecture,
see~\cite{BK03a,BK03b,BK03c,IAS}. Since then the study of $\thom$
complexes in general has become a~vibrant field, see additionally
\cite{Cs,CK1,CK2,En05,Ka05,K4,K5,K06,Pf05,Sch,Sch2}. The $\thom$
complexes even appeared in topological stochastics, see
\cite{MSS03,SS}.

The particular information about the cohomology groups of
$\thom(C_m,K_n)$, which was needed for the proof of this conjecture
consisted of computing the cohomology groups, with $\zz$ as well as
with integer coefficients, up to dimension $n-2$, analyzing the
$\zz$-action on the cohomology groups (in the case of
$\dz$-coefficients), and computing the height of the Stiefel-Whitney
characteristic classes (for odd~$n$).

Here we present a~complete computation of the groups $H^*(C_m,K_n)$
first with $\zz$, and then with integer coefficients. The idea is to
take the same spectral sequence as in \cite{BK03c}, but this time to
calculate the second tableau in its entirety. The main technical tool
will be the detailed analysis of a~certain family of combinatorially
defined cubical complexes $\{\Phi_{m,n,g}\}$, which is interesting in
its own right. Once these complexes are understood, we will almost
have computed all the cohomology groups. The two additional
ingredients are the computation of signs, to handle the case of the
integer coefficients, as well as in-depth study of the top row entries
in the second tableau.

There are two byproducts of our investigation. First, several
technical lemmas from \cite{BK03c} are special cases of our statements
here. Thus, we find alternative simpler ways to confirm their
correctness. Second, the collapsing procedure for the complexes
$\Phi_{m,n,g}$, which we call {\it grinding} for the graphic reasons
which will become clear later, generalizes easily to yield
a~collapsing procedure for the connected components of complexes
$\thom(C_m,C_n)$, leading to symmetrical cores consisting of garlands
of higher-dimensional cubes, linked in a~circle pattern by their
opposite end vertices. Thereby we obtain a~combinatorial strengthening
of the earlier result that these connected components are either
points or are homotopy equivalent to circles, see~\cite{CK1}.

\vskip5pt

\nin The paper is organized as follows.

\vskip5pt

\nin {\bf Section 2.} This section contains preliminaries on $\thom$ 
and $\thomp$ complexes of arbitrary graphs, as well as on independence
and partial coloring complexes of cycles.

\vskip5pt

\nin {\bf Section 3.} Here, we review the setup of the spectral 
sequence converging to the cohomology of $\thomp(T,G)$, and recall
various facts about the first two tableaux, and the first differential
of that sequence.

\vskip5pt

\nin {\bf Section 4.} We introduce the terminology of arc pictures 
and see how the generators of the entries in the first tableau of our
spectral sequence can be indexed by them. This allows us to perform
the first reduction in our computation of the second differential.

\vskip5pt

\nin {\bf Section 5.} The bulk of our calculation is concentrated in 
this section and in Section~\ref{s:zsect}. We rephrase the problem of
computing the entries of the second tableau as a~question of computing
the homology groups of a~certain family of cubical complexes, the
so-called {\it torus front complexes}. We find a~collapsing procedure
on these complexes, which we call {\it grinding}, as a~result reducing
them to the thin torus front complexes, which turn out to be garlands
of cubes linked by the opposite end vertices. This allows almost
a~complete computation of the second tableau of our spectral sequence
with $\zz$-coefficients.

\vskip5pt

\nin {\bf Section 6.} We use our filtration to prove a formula which 
expresses the Euler characteristic of the $\thom$ complex in terms of
Euler characteristics of $\thomp$ complexes. Specifically, the Euler
characteristic of $\thom(T,K_n)$ can be expressed in terms of Euler
characteristics of the independence complexes of induced subgraphs
of~$T$. This leads to an easy proof of an older formula. Furthermore,
we use the second tableau of our spectral sequence to derive a~formula
for the Euler characteristic of $\thom(C_m,K_n)$.

\vskip5pt

\nin {\bf Section~\ref{s:zsect}.} The complete calculation of the groups
$H^*(\thom(C_m,K_n);\dz)$ is finalized. Since here we are working over
integers, the additional element of computing signs in the
differentials requires some care, as does the analysis of a~few top
diagonals in the spectral sequence tableaux.

\vskip5pt

\nin {\bf Section 8.} We apply the grinding technique to demonstrate
a~useful collapsing of $\thom(C_m,C_n)$ onto its symmetric core of
thin cycle maps, which also turns our to be a~collection of garlands
of cubes.


\section{Preliminaries}

\subsection{$\thom$ versus $\thomp$.} \label{ssect:prel1} $\,$

\vskip5pt

\nin All graphs considered in this paper will be finite.  Let us start with
the basic definitions. Let $G$ be a~graph, and let $A,B\subseteq
V(G)$, $A,B\neq\emptyset$. We call $(A,B)$ a~{\it complete bipartite
subgraph} of $G$, if for any $x\in A$, $y\in B$, we have $(x,y)\in
E(G)$, i.e., $A\times B\subseteq E(G)$.

\begin{df} \label{df:hom}
For arbitrary graphs $T$ and $G$, the prodsimplicial complex
$\thom(T,G)$ is the subcomplex of $\prod_{x\in V(T)} \Delta^{V(G)}$
defined by the following condition: $\sigma=\prod_{x\in
V(T)}\sigma_x\in\thom(T,G)$ if and only if for any $x,y\in V(T)$, if
$(x,y)\in E(T)$, then $(\sigma_x,\sigma_y)$ is a~complete bipartite
subgraph of~$G$.
\end{df}

Here, $\prod_{x\in V(T)} \Delta^{V(G)}$ denotes the~direct
product of $|V(T)|$ identical copies of $\Delta^{V(G)}$, indexed by
the vertices of~$T$.  We refer the reader to the survey article
\cite{IAS} for all further background material on $\thom$ complexes.

Allowing partial graph homomorphisms yields a~related simplicial
complex.

\begin{df} \label{df:homp}
For arbitrary graphs $T$ and $G$, $\thomp(T,G)$ is the simplicial
subcomplex of $\join{x\in V(T)} \Delta^{V(G)}$ defined by the
following condition: $\sigma=\join{x\in V(T)}\sigma_x\in\thomp(T,G)$
if and only if for any $x,y\in V(T)$, if $(x,y)\in E(T)$, and both
$\sigma_x$ and $\sigma_y$ are nonempty, then $(\sigma_x,\sigma_y)$ is
a~complete bipartite subgraph of~$G$.
\end{df}

We call $\thomp(T,G)$ the {\it Hom plus complex}. The plus version of
the space of all colorings of a~certain graph has a~simpler
topological structure, which we now proceed to describe.

For an arbitrary graph $G$, let the {\it strong complement} $\cmp G$
be the graph defined by $V(\cmp G)=V(G)$, and $E(\cmp G)=V(G)\times
V(G)\sm E(G)$. Let the {\it independence complex} $\ind(G)$ be the
simplicial complex whose set of vertices is $V(G)$, defined by the
rule: $\sigma\in 2^{V(G)}$ is a~simplex in $\ind(G)$ if and only if
$\sigma$ is an~independent set of $G$, i.e., no two vertices in
$\sigma$ are connected by an~edge.

\begin{prop} \label{pr:homp}
{\rm (\cite[Proposition 3.2]{BK03c}).}

\nin For arbitrary graphs $T$ and $G$,
$\thomp(T,G)$ is isomorphic to $\ind(T\times\cmp G)$. 
\end{prop}

Specializing Proposition~\ref{pr:homp} to the case $G=K_n$, and taking
into account that $\cmp K_n= \coprod_{i=1}^n \cmp K_1$, as well as the
fact that for arbitrary graphs $G_1$ and $G_2$ we have the isomorphism
of simplicial complexes
\[\ind\left(G_1\coprod G_2\right)=\ind(G_1)*\ind(G_2),
\]
 we obtain the following corollary.

\begin{crl} \label{crl:hompkn}
For an~arbitrary graph $T$,
$\thomp(T,K_n)$ is isomorphic to the $n$-fold join $\ind(T)^{*n}$.
\end{crl}

In other words, if we understand the independence complex of a~graph,
then we understand the plus version of its space of $n$-colorings for
any~$n$.

An important additional piece of structure, binding together the Hom
and the Hom plus constructions is that of a~support map. Indeed, for
any topological space $X$ and a set $S$, there is the standard support
map from the $S$-fold join of $X$ to the appropriate simplex $\supp:
*_S X\longrightarrow \Delta^S$.  The map $\supp$ is simply the
projection map which "forgets" the coordinates in $X$.

As a~special case, for arbitrary graphs $T$ and $G$, we get the
restriction map $\supp:\thomp(T,G)\ra\Delta^{V(T)}$. Explicitly, for
each simplex of $\thomp(T,G)$, $\eta:V(T)\ra 2^{V(G)}$, the support of
$\eta$ is given by $\supp\eta=V(T)\setminus\eta^{-1}(\emptyset)$.
An~important property of the support map is that the preimage of the
barycenter of $\Delta^{V(T)}$ is homeomorphic to $\thom(T,G)$. We
shall use this heavily in the analysis of our spectral sequence.

\subsection{Cohomology of partial colorings of cycles.} \label{ssect:prel2} $\,$

\vskip5pt

\nin The homotopy type of the independence complexes of cycles is well-understood.

\begin{prop} \label{pr:indcr}
{\rm (\cite[Proposition 5.2]{K2})}.

\nin For any $m\geq 2$, we have
\[\ind(C_m)\simeq\begin{cases}
S^{k-1}\vee S^{k-1},&\text{ if } m=3k;\\
S^{k-1},&\text{ if } m=3k\pm 1.
\end{cases}\]
\end{prop}
By Corollary~\ref{crl:hompkn}, we have $\thomp(C_m,K_n)\simeq\ind(C_m)^{*n}$, 
hence we can derive an~explicit description.


\begin{crl} \label{crl:cplus}
{\rm (\cite[Corollary 4.2]{BK03c})}.

\nin For any $m\geq 2$, $n\geq 3$, we have
\nin 
        \[\thomp(C_m,K_n)\simeq
        \begin{cases}
        \bigvee_{2^n\text{ copies}} S^{nk-1},&\text{ if } m=3k;\\
        S^{nk-1},&\text{ if } m=3k\pm 1.
        \end{cases}\]
\end{crl}

We shall use this result later on, as our spectral sequence will be
converging to the cohomology groups of $\thomp(C_m,K_n)$. For future reference
we also record the following observation.

\begin{prop} \label{pr:mnineq}
For all integers $m,n$, such that $m\geq 5$, $n\geq 4$, we have $\wti
H^i(\thomp(C_m,K_n))=0$, for $i\leq m+n-4$, with one exception $(m,n,i)=(7,4,7)$.
\end{prop}

\pr By Corollary \ref{crl:cplus}, we just need to check that 
\begin{equation} \label{eq:m+n}
 nk-1>m+n-4, 
\end{equation}
where $k=\lfloor(m+1)/3\rfloor$.
  
If $m\in\{5,6\}$, then $k=2$, and \eqref{eq:m+n} reduces to $n+3>m$,
which is true.  Thus we can assume that $m\geq 7$. Notice that
$k\geq(m-1)/3$, hence \eqref{eq:m+n} would follow from
$n(m-1)/3+3>m+n$. This can be rewritten as $nm-3m-4n+9>0$, or,
equivalently as $(n-3)(m-4)>3$. Since $n\geq 4$, $m\geq 7$, this
inequality is valid, with the equality only when $m=7$, $n=4$.
\qed

\section{Setting up the spectral sequence}

\subsection{Filtration induced by the support map.} \label{ssect:supp} $\,$ \vskip5pt


\nin For simplicity, we shall for now consider $\zz$-coefficients. The additional issues arising in connection with switching to integer coefficients will be analyzed in Section~\ref{s:zsect}. 

For a~function $\eta:V(T)\ra 2^{V(G)}$ let $\eta_+$ denote the chain in $C_*(\thomp(T,G))$ consisting of a~single simplex (with coefficient 1) indexed by $\eta$; when no confusion arises, we identify this chain with the simplex itself. Furthermore, let $\eta^*_+$ denote the dual cochain.
Following \cite{BK03c,IAS} we consider the Serre filtration of the cellular cochain complex $C^*(\thomp(T,G);\zz)$ associated with the support map. To describe the considered filtration explicitly,  define the subcomplexes $F^p=F^pC^*(\thomp(T,G);\zz)$ of $C^*(\thomp(T,G);\zz)$ as follows:
\[ F^p: \dots\stackrel{\bo^{q-1}}\lra F^{p,q}\stackrel{\bo^{q}}\lra
F^{p,q+1}\stackrel{\bo^{q+1}}\lra\dots,\]
where 
\[ F^{p,q}=F^p C^q(\thomp(T,G);\zz)=\zz\left[\eta^*_+\,\left|\,
\eta_+\in\thomp^{(q)}(T,G),|\supp\eta|\geq p+1\right]\right.,
\]
$\bo^*$ is the restriction of the differential in
$C^*(\thomp(T,G);\zz)$, and $\thomp^{(q)}(T,G)$ denotes the $q$-th
skeleton of $\thomp(T,G)$. The idea is to filter by the preimages of
the skeleta of $\Delta^{V(T)}$. The $F^{p,q}$ is a~vector space over $\zz$
generated by all elementary cochains corresponding to
$q$-dimensional cells, which are supported in at least $p+1$ vertices
of~$T$. Note, that this restriction defines a~filtration, since the
differential does not decrease the cardinality of the support set.
 We have
\[
C^q(\thomp(T,G);\zz)=F^{0,q}\supseteq F^{1,q}\supseteq\dots\supseteq
F^{|V(T)|-1,q}\supseteq F^{|V(T)|,q}=0.
\]

\subsection{The 0th and the 1st tableaux.} \label{ssect:01tabl} $\,$ \vskip5pt


\nin We follow the treatment in \cite{BK03c,IAS} to describe the 0th and the 1st tableaux, as well as to perform a~partial analysis of the second tableau.
 
To start with, as an~additional piece of notations, by writing the brackets $[-]$ after the name of a~cochain complex, we shall mean the index shifting (to the left), that is for the cochain complex $\cc=(C^*,d^*)$, the cochain complex $\cc[s]=(C^*[s],d^*)$ is defined by $C^i[s]:=C^{i+s}$. We sharpen reader's attention on the fact that, when considering integer coefficients, we choose not to change the sign of the differential.

\begin{prop} {\rm(\cite[Proposition 3.4]{BK03c}).}

\nin For any $p$, we have
\begin{equation}\label{eq:e0tab}
F^p/F^{p+1}=\bigoplus_{\begin{subarray}{c}{S\subseteq V(T)}\\{|S|=p+1}
\end{subarray}}
C^*(\thom(T[S],G);\car)[-p].
\end{equation}
Hence, the 0th tableau of the spectral sequence associated to the
cochain complex filtration $F^*$, and converging to
$H^{p+q}(\thomp(T,G);\car)$, is given by
\begin{equation}
  \label{eq:E0}
E_0^{p,q}=C^{p+q}(F^p,F^{p+1})=\bigoplus_{\begin{subarray}{c}
{S\subseteq V(T)}\\{|S|=p+1}\end{subarray}}
C^q(\thom(T[S],G);\car).
\end{equation}
\end{prop}

Furthermore, using the standard facts about spectral sequence, see~\cite{IAS,McC}, we obtain the description of the first tableau as well.

\begin{equation} \label{eq:E1}
E_1^{p,q}=H^{p+q}(F^p,F^{p+1})=\bigoplus_{\begin{subarray}{c}
{S\subseteq V(T)}\\{|S|=p+1}\end{subarray}}H^q(\thom(T[S],G);\car).
\end{equation}

\subsection{The first differential.} \label{ssect:1diff} $\,$ \vskip5pt

 
\nin  According to the formula \eqref{eq:E1} the first differential $d_1^{p,q}:E_1^{p,q}\lra E_1^{p+1,q}$ is actually a~map
\[\bigoplus_{\begin{subarray}{c}{S\subseteq V(T)}\\{|S|=p+1}\end{subarray}} H^q(\thom(T[S],G);\zz)\lra \bigoplus_{\begin{subarray}{c}
{S\subseteq V(T)}\\{|S|=p+2}\end{subarray}}H^q(\thom(T[S],G);\zz).
\]

For $S_2\subseteq S_1\subseteq V(T)$, let $i[S_1,S_2]:T[S_2] \hookrightarrow T[S_1]$ be the inclusion graph homomorphism. Since $\thom(-,G)$ is a~contravariant functor, we have an~induced map $i_G[S_1,S_2]:\thom(T[S_1],G)\ra\thom(T[S_2],G)$, and hence, an~induced map on the cohomology groups 
\[
i_G^*[S_1,S_2]:H^*(\thom(T[S_2],G);\car)\ra H^*(\thom(T[S_1],G);\car).
\]
Let $\sigma\in H^q(\thom(T[S],G);\zz)$, for some $q$, and some $S\subseteq V(G)$. The value of the first differential on $\sigma$ is then given by

\begin{equation} \label{eq:d1pq}
d_1^{p,q}(\sigma)=\sum_{x\in V(T)\sm S} i_G^*[S\cup\{x\},S](\sigma).
\end{equation}

We shall next detail the formula \eqref{eq:d1pq} by introducing specific generators for the groups in the first tableau.

\section{Encoding cohomology generators by arc pictures}

\subsection{The language of arcs.} \label{ssect:arc1} $\,$ \vskip5pt

\nin An extensive notational apparatus to deal with these specific 
spectral sequence computations was introduced in \cite{BK03c}. It is
convenient here to alter these notations somewhat, in part since our
needs are different and in part to streamline the presentation, so,
for the sake of compatibility, we shall now introduce the notations
which we need.

Let $S$ be a~proper subset of $V(C_m)$. We call the connected
components of the graph $C_m[S]$ either {\it singletons} or {\it arcs}
depending on whether they have 1 or at least 2 vertices. For $v\in S$
we let $a(S,v)$ denote the arc of $S$ to which $v$ belongs (assuming
this arc exists). Furthermore, for an~arbitrary arc $a$ of $S$, we let
$a=[a_\bu,a^\bu]_m$, and finally, we set $\widehat
a:=[a_\bu-1,a^\bu+1]_m$, so $|\widehat a|=|a|+2$, if $|a|\leq m-2$. We
stress that the arithmetic operations as well as intervals are taken
modulo $m$, for example $[m,1]_m$ and $[1,2]_m$ are arcs with two
vertices each, while $\widehat{[m,1]}_m=[m-1,2]_m$ and
$\widehat{[1,2]}_m=[m,3]_m$. To make the formulas easier, we also
think of the cycle $a=C_m$ itself as an arc, in which case we use the
convention $a_\bu=1$. When $A$ is a~collection of arcs, we set
$V(A):=\bigcup_{a\in A}a$, $A_\bu:=\bigcup_{a\in A}a_\bu$, and we let
$\widehat A$ denote the union $\bigcup_{a\in A}\widehat a$.

\begin{df}
Let $t\geq 1$. A~{\bf $t$-arc picture} is a~pair $(S,A)$, where $S$ is some proper subset of $V(C_m)$, and $A$ is a~set consisting of $t$ different arcs of $C_m[S]$. We refer to arcs in $A$ as {\bf marked arcs}, and sometimes simply call $(S,A)$ an~{\bf arc picture}.
\end{df}

\subsection{The corresponding cohomology generators.} \label{ssect:arc2} $\,$ \vskip5pt

\nin For $S\subset V(C_m)$, such that $|S|=p+1<m$, the $t$-arc pictures $(S,A)$ will index the generators of $E_1^{p,t(n-2)}$. To introduce these, let $V=\{v_a\}_{a\in A}$ be a~set of vertex representatives of $A$, i.e., $v_a\in A$ for each $a\in A$. Now, set
\begin{equation} \label{eq:eta+}
\sigma_V^S:=\sum_\eta \eta^*_+,
\end{equation}
where the sum is taken over all $\eta:V(S)\ra 2^{[n]}\sm\{\emptyset\}$, such that
$\eta(v_a)=[n-1]$, for all $a\in A$, and $|\eta(w)|=1$, for all $w\in V(S)\sm V$.

It is easily checked that $\sigma_V^S$ is a~cocycle of
$\thomp(C_m[S],K_n)$.  When vertex representatives move along edges,
the corresponding cohomology class $\left[\sigma_V^S\right]$ does not
change, see
\cite[subsection~4.2]{BK03c}, therefore, we may limit our attention to
$\sigma_{A_\bu}^S$. The only reason we introduce the vertex
representatives at all is for the later calculations with integer
coefficients. In this case, the cohomology class $\left[\sigma_V^S\right]$
changes the sign according to a~certain pattern, see
Section~\ref{s:zsect} of this paper and \cite[subsection~4.2]{BK03c}.

It is now time to interpret the formula \eqref{eq:d1pq} in terms of
our combinatorial generators. Since the differentials in the spectral
sequence are induced by the differential in the original chain
complex, we may conclude that
\begin{equation} \label{eq:d1gen}
d_1(\left[\sigma_{A_\bu}^S\right])=\sum_{w\notin S}\left[\sigma_{A_\bu}^{S\cup\{w\}}\right],
\end{equation}
where the sum is taken over all $w$, such that adding $w$ to $S$ does not unite two arcs from $A$. Note that \eqref{eq:d1gen} is also valid when $|S|=m-1$, if we define $\sigma_V^S$ in the same way for $S=[m]$. It might be worthwhile to observe at this point that the set of marked arcs in $S\cup\{w\}$ is equal to $A$ if and only if $w\notin\widehat{A}$.

\subsection{The first reduction.} \label{ssect:arc3} $\,$ \vskip5pt

\nin By the definition (or in some approaches - property) of spectral sequences, in order to obtain the second tableau $\{E_2^{*,*}\}$ we need to calculate the cohomology groups of the complex $(E_1^{*,(n-2)t},d_1)$. Let $(A_t^*,d_1)$ be its subcomplex generated by all classes $\left[\sigma^S_{A_\bu}\right]$, such that $\widehat{A}=V(C_m)$.

\begin{prop} \label{pr:red1} {\rm (\cite[Lemma 4.8]{BK03c}).} $\,$

\nin For all $t$, we have
\begin{equation} \label{eq:red1}
 H^*(E_1^{*,(n-2)t})=H^*(A_t^*).
\end{equation}
\end{prop}

To be precise, Proposition \ref{pr:red1} is only stated in
\cite[Lemma~4.8]{BK03c} for the case when $m$ is odd. However, the
parity of $m$ is not used in the proof there, and therefore holds for
all~$m$. The proof of \cite[Lemma~4.8]{BK03c} proceeds by setting up
an~auxiliary spectral sequence and then analyzing it. It is also
possible to give a~more direct combinatorial (matching) argument by
using the version of the discrete Morse theory which works for chain
complexes, see~\cite{dmt}. 

The matchings can be done as follows: for each collection of arcs $A$,
such that $\widehat{A}\neq V(C_m)$, choose some element $x\in
V(C_m)\sm\widehat A$, then, an~element $\left[\sigma_{A_\bu}^S\right]$
is matched with the element
$\left[\sigma_{A_\bu}^{S\xor\{x\}}\right]$, where $\xor$ denotes the
symmetric difference (exclusive or) of sets. One can then check that
this matching is acyclic, in the sense of~\cite{dmt}. Indeed, assume
that there exists a~cycle $y_t\prec x_1\succ y_1\prec x_2\succ\dots
\prec x_t\succ y_t$, where $x_i$ and $y_i$ are matched for all 
$i=1,\dots,t$. Consider the covering relation $x_i\succ y_{i-1}$, for some
$i=2,\dots,t$. One can see that in this situation, $x_i$ must be
obtained from $y_{i-1}$ by adding an element to its indexing set $S$
which extends one of the arcs in $A$. The same should hold for the
covering relation $y_t\prec x_1$. Clearly, this leads to a~contradiction, since the total length of arcs in $A$ does not change on the matching edges
in the cycle ($A$ itself does not change there), while it increases
along other edges.

\section{Topology of the torus front complexes} \label{sect5}

\subsection{Reinterpretation of $H^*(A_t^*,d_1)$ using a~family 
of cubical complexes $\{\Phi_{m,n,g}\}$.} \label{ssect:tor1} $\,$
\vskip5pt

\nin Combining Proposition \ref{pr:red1} with formula \eqref{eq:d1gen} 
we see that to compute the second tableau of our spectral sequence we
need to calculate cohomology groups of a~certain cochain complex,
which we have combinatorially described by choosing an~explicit basis of
generators and then writing out the differentials in terms of these
generators.

Our next idea is to reinterpret the cochain complex $(A_t^*,d_1)$ as
a~chain complex which computes $\zz$-homology of a~certain cubical
complex. We now proceed with defining these complexes.  For a~natural
number $n\geq 2$, a~circular $n$-set is the set of $n$ points, which
are equidistantly arranged on the circle of length~$n$.

\begin{df} \label{df:phi}
Let $m,n$, and $g$ be natural numbers. The cubical complex
$\Phi_{m,n,g}$ is defined as follows:
\begin{itemize}
\item vertices of $\Phi_{m,n,g}$ are indexed by all possible 
selections of $m$ elements from the circular $n$-set, so that every
two of the chosen vertices are at a~distance at least~$g$;
\item to index higher-dimensional cubes we take all possible 
collections of $m$ sets, where each set either contains a~single
element of the circular $n$-set, or consists of a~pair of two elements
at distance~1; the sets are requested to be at minimal distance $g$,
where the distance between two sets is defined as the minimum of the
distances between their elements.
\end{itemize}
\end{df}

\begin{figure}[hbt]
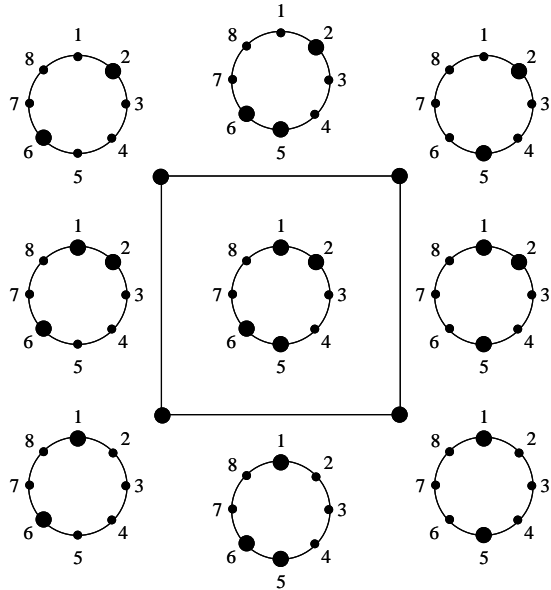

\begin{center}
  \begin{picture}(0,0)%
    \includegraphics{phidef.pstex}%
  \end{picture}%
  \input{phidef.pstex_t}%
  
\end{center}
\caption{Indexing of cubes in $\Phi_{2,8,2}$.}
\label{fig:phidef}
\end{figure}

Clearly, the dimension of a~cube indexed by a~collection of sets is
equal to the number of 2-element sets in this collection, and one cube
contains another if and only if the corresponding collection of sets
contains the collection of sets associated to the second cube. See
Figure~\ref{fig:phidef} for graphic explanation.

For convenience of formulations, we introduce the convention that the
complex $\Phi_{0,n,g}$ is a~point. One can quickly list some trivial
cases of the complexes $\Phi_{m,n,g}$.  To start with a~necessary and
sufficient condition that these complexes are non-empty is that $n\geq
gm$. The complex $\Phi_{n,n,1}$ is just a~point, and more generally
$\Phi_{m,gm,g}$ is a~disjoint union of $g$ points. It is not difficult
to check by hand that the complex $\Phi_{m,n,g}$ is connected whenever
$n>gm$. The next case to consider is $\Phi_{m,gm+1,g}$, the interested
reader is invited to check that $\Phi_{m,gm+1,g}$ is a~cycle of length
$gm+1$. Finally we remark that the dimension of $\Phi_{m,n,g}$ is
equal to $\min(m,n-gm)$.

The usefulness of this family in our context becomes apparent from the
following proposition.

\begin{prop} \label{pr:tm2}
The cochain complex $(C_*(\Phi_{t,m,3};\zz),d)$ is isomorphic to the
chain complex $(A_t^*,d_1)$, with module $C_i(\Phi_{t,m,3};\zz)$
corresponding to the module $A_t^{m-t-i-1}$, for all~$i$.
\end{prop}

\pr We think of a~collection of sets which indexes a~cube in 
$\Phi_{t,m,3}$ as a~collection of gaps between the arcs from $A$. This
will clearly give a~dimension-preserving bijection between the cubes
of $\Phi_{t,m,3}$ and the collections of $t$ arcs $A$ such that
$\widehat A=V(C_n)$. The condition that the distance between sets
(read - gaps) is at least~3 corresponds to the requirement that arcs
must have at least~2 elements. It remains the compare the differential
maps: they coincide by the formula~\eqref{eq:d1gen}. Note that we are
using the fact that we are working over $\zz$-coefficients.
\qed

\vskip5pt

Since we know that $\Phi_{m,gm,g}$ are disjoint unions of $g$ points
it is enough to assume from now on that $n>gm$.

\subsection{The torus front interpretation.} \label{ssect:tor2} $\,$ \vskip5pt

\nin Let $m$ and $n$ be two positive integers. An {\it $(m,n)$-grid path}
is a~directed path on the unit orthogonal grid using only the unit
basis vectors $(0,1)$ (northbound edge) and $(1,0)$ (eastbound edge)
and connecting the point $(-x,x)$ with the point $(m-x,n+x)$, for some
integer~$x$.

\begin{df} \label{df:torfr}
Consider the action of the group $\dz\times\dz$ on the plane, where
the standard generators of $\dz\times\dz$ act by vector translations,
with vectors $(-1,1)$ and $(m,n)$. An {\bf $(m,n)$-torus front} (or
sometimes simply torus front) is an~orbit of this
$(\dz\times\dz)$-action on the set of all $(m,n)$-grid paths.
\end{df}

\begin{figure}[hbt]
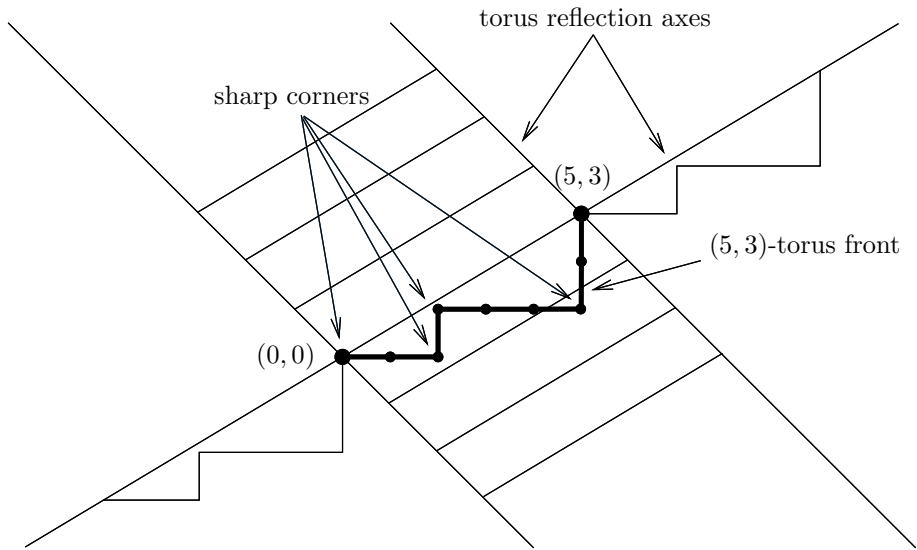

\begin{center}
  \begin{picture}(0,0)%
    \includegraphics{torfr.pstex}%
  \end{picture}%
  \input{torfr.pstex_t}%
  
\end{center}
\caption{$(5,3)$-path and the associated torus front.}
\label{fig:torfr}
\end{figure}

The reason for choosing the word {\it torus} here is because the
quotient space of the plane by this $(\dz\times\dz)$-action can be
viewed as a~torus, where each $(m,n)$-grid path yields a~loop
following the grid in the northeastern direction, see
Figure~\ref{fig:torfr}. We notice that every torus front has a~unique
representative starting from the point~$(0,0)$. It will soon become
clear why we choose to consider the orbits of the action rather than
merely considering these representatives of the orbits.

A {\it northwestern sharp corner} of a~torus front is a~vertex which
is entered by a~northbound edge, and exited by an~eastbound edge, in
the same way, a~{\it southeastern sharp corner} of a~torus front is
a~vertex which is entered by an~eastbound edge, and exited by
a~northbound edge.

An {\it elementary flip} of a~torus front is either a~replacement of
a~northwestern sharp corner at $(a,b)$ by a~southeastern sharp corner
at $(a+1,b-1)$ or vice versa. It may help to think a~torus front as
a~sort of flexible snake, where the sharp corners (meaning the vertex
in the sharp corner, together with the two adjacent edges) can be
flipped about (as if using ball-and-socket joints) the diagonal line
connecting the neighbors of the sharp corner vertex. 

\begin{figure}[hbt]
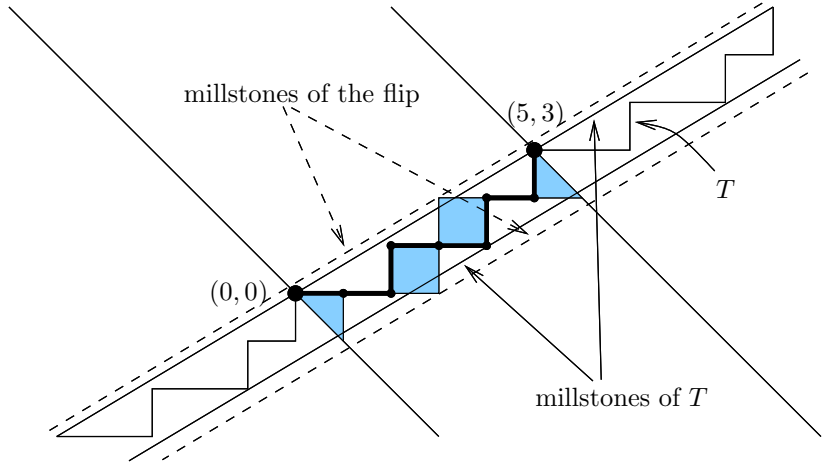

\begin{center}
  \begin{picture}(0,0)%
    \includegraphics{flip1.pstex}%
  \end{picture}%
  \input{flip1.pstex_t}%
  
\end{center}
\caption{A flip of a $(5,3)$-torus front consisting of 3 non-interfering elementary flips.}
\label{fig:flip1}
\end{figure}

Finally, A~{\it flip} of a~torus front is a~collection of
non-interfering elementary flips, i.e., the flipped vertices are at
least at distance~2 from each other, where the distance is measured along the front, see
Figure~\ref{fig:flip1}. The picture which we have in mind is that the torus front propagates through the torus by means of flips. Since we want flips to encode cells, we prefer to think of them as schemes of allowed elementary flips, rather than a~concrete process of replacing one torus front with another. For a~flip $F$, we call those torus fronts, which can be obtained by actually performing all the elementary flips from $F$ (meaning here choosing one of the two positions for each flipped corner), the {\it members} of~$F$.

\begin{df} \label{df:tfcomp}
Let $m,n$, and $g$ be natural numbers. The cubical complex
$TF_{m,n,g}$ is defined as follows:
\begin{itemize}
\item vertices of $TF_{m,n,g}$ are indexed by all possible 
$(m,n)$-torus fronts, whose horizontal legs have length at least~$g$, here the length of the leg is taken to be the number of vertices it contains;
\item the higher-dimensional cubes of $TF_{m,n,g}$ are indexed by all 
possible flips of $(m,n)$-torus fronts, whose members are vertices of $TF_{m,n,g}$.
\end{itemize}
\end{df}

Notice that for $g=1$ there are no conditions, so $TF_{m,n,1}$ simply
has all $(m,n)$-torus fronts as vertices and all flips as
higher-dimensional cubes. It is also important to remark that the term "horizontal legs" includes the legs of length $1$, on other words, for $g>1$, it is prohibited that the torus front contains two consecutive northbound edges.

\begin{prop}\label{pr:tfint}
For any natural numbers $m,n$, and $g$, the cubical complexes
$\Phi_{m,n,g}$ and $TF_{m,n-m,g}$ are isomorphic.
\end{prop}

\pr Let $S\subseteq [n]$, $|S|=m$, be a~vertex of $\Phi_{m,n,g}$. 
We can construct an~$(m,n-m)$-grid path as follows. Start out from the
point $(0,0)$, then, for all $i$ from $1$ to $n$ perform the following
step: either move along the northbound edge if $i\in S$, or move along
the eastbound edge if $i\notin S$. Since $|S|=m$, we will eventually
reach the point $(m,n-m)$. It is obvious that this gives a~bijection
between vertices of $\Phi_{m,n,g}$ and $(m,n-m)$-grid paths satisfying
the extra restriction that the horizontal legs are at least $g$-long.

It is also clear that moving along edges in $\Phi_{m,n,g}$ corresponds
to elementary flips in $TF_{m,n-m,g}$. This is completely obvious as
long as the flip does not concern the vertex $(0,0)$ (or, equivalently
the vertex $(m,n-m)$). If, on the other side, the edge in
$\Phi_{m,n,g}$ indicates that we should flip the sharp corner $(0,0)$,
then we can do that as well, but then we should also flip $(m,n-m)$,
and we will get an~$(m,n-m)$-path from $(1,-1)$ to $(m+1,n-m-1)$;
translating back we will again get an~$(m,n-m)$-path from $(0,0)$ to
$(m,n-m)$.

If the flips are done on torus fronts instead of grid paths, then we
do not need to consider different special cases, which is the main
reason why we chose to replace the paths by the orbits of the
$(\dz\times\dz)$-action.
\qed

\vskip5pt

Note that the complex $TF_{m,n-m,1}$ is isomorphic to the complex
$TF_{n-m,m,1}$, and so the complex $\Phi_{m,n,1}$ is isomorphic to the
complex $\Phi_{n-m,n,1}$. 

To understand the complexes $\Phi_{m,n,g}$ satisfying $n>mg$, we shall
study the complexes $TF_{m,n,g}$ with $n>m(g-1)$.

\subsection{Grinding.} \label{ssect:tor3} $\,$ \vskip5pt

\begin{figure}[hbt]
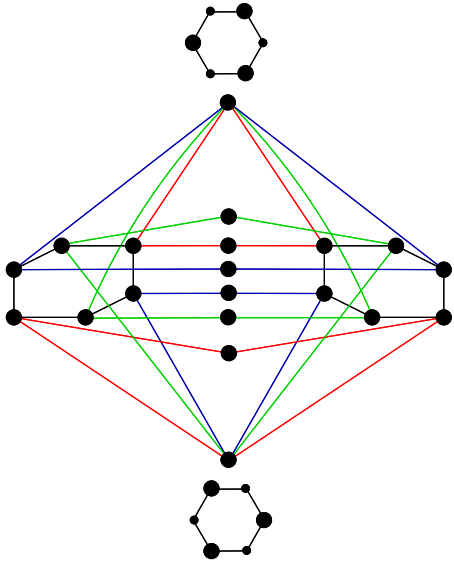

\begin{center}
  \begin{picture}(0,0)%
    \includegraphics{phi361.pstex}%
  \end{picture}%
  \input{phi361.pstex_t}%
  
\end{center}
\caption{The complex $TF_{3,3,1}=\Phi_{3,6,1}$.}
\label{fig:phi361}
\end{figure}

\nin We shall now drastically simplify complexes $TF_{m,n,g}$ be means of 
the collapsing procedure which we call {\it grinding}. To illustrate
the idea, let us consider the case $m=n=3$, $g=1$. The complex
$TF_{3,3,1}$ ($=\Phi_{3,6,1}$) consists of two cubes, indexed by the
set collections $\{\{1,2\},\{3,4\},\{5,6\}\}$ and
$\{\{1,6\},\{2,3\},\{4,5\}\}$, and six additional squares, see
Figure~\ref{fig:phi361}\footnote{The figure is reproduced courtesy to
Eva-Maria Feichtner.}. Each of these squares commands a~vertex which
does not belong to any other maximal cell; for example, for the square
indexed with $\{\{1,2\},\{3\},\{4,5\}\}$ this vertex is indexed with
$\{\{2\},\{3\},\{4\}\}$. Therefore, we can collapse these 6 squares
through these 6 vertices and the only thing left will be the 2 cubes
hanging together by their end points. This is the characteristic
picture we shall now see in general.

Let $T$ be an $(m,n)$-torus front, which we represent as an
$(m,n)$-grid path starting from $(0,0)$. We associate to $T$ two lines
$M_{nw}(T)$ and $M_{se}(T)$ defined as follows: both lines are
parallel to the line $L$ passing through the points $(0,0)$ and
$(m,n)$, and they both are tangential to the grid path representing
$T$, $M_{nw}(T)$ from the northwest, and $M_{se}(T)$ from the
southeast; see Figure~\ref{fig:flip1}. We call these two lines the
{\it millstones} of~$T$.

The {\it width} of the torus front $T$, denoted $w(T)$, is the distance between the millstones measured along the line $x=-y$. Clearly, the width cannot be less than $\sqrt{2}/2$, and the minimum is achieved only when $m=n$ by the torus front where northbound and eastbound edges alternate at every step.

We let $E_{nw}(T)$ denote the set of the points in $T\cap M_{nw}(T)$ which we call {\it northwestern extreme corners}; in the same way $E_{se}(T)=T\cap M_{se}(T)$ denotes the set of southeastern extreme corners. Extreme corners are by construction always sharp corners.

As mentioned above, each sharp corner defines an elementary flip of $T$, during which it remains a~sharp corner, but changes orientation. If this sharp corner was extreme, then the new sharp corner would be extreme (on the opposite side) if and only if $w(T)\leq\sqrt 2$. Clearly, for any torus front, the elementary flips of northwestern extreme corners are mutually non-interfering, and therefore form one flip; the same is true for the southeastern extreme corners. Furthermore, if $w(T)>\sqrt2/2$, then these two sets of elementary flips do not interfere with each other either. Indeed, if they did, it would mean that two neighboring vertices on the torus front are extreme corners of different orientations, which would imply $w(T)=\sqrt2/2$.

Assume now that $w(T)>\sqrt 2$, and consider the total flip of all extreme corners. It is geometrically clear, that performing any combination of these elementary flips will not increase the width of the torus front, and that performing all of the northwestern flips, or all of the southeastern flips, simultaneously will for sure decrease it. 

Another important observation is that if an extreme northwestern corner is flipped on a torus~front, where all horizontal legs have length at least~$g$, then all horizontal legs have length at least~$g$ in the obtained torus front as well; the same is true if we flip an extreme southeastern corner. To see this, note that the only time the horizontal leg could be decreased, is when one of its endpoints would be flipped. Assume that the horizontal leg in question runs from vertex $a$ to vertex $b$, that $a$ is flipped (no loss of generality here), and denote by $c$ the vertex which we reach by the northbound edge from~$b$. Then, the slope of the line connecting $a$ with $c$ is $1/(g-1)$, and, since we assumed $m(g-1)<n$, we have $1/(g-1)>m/n$. The vertex $a$ is flipped, hence it is a~northwestern extreme corner, since $c$ must lie on $M_{nw}(T)$ or to the southeast of it, this implies that $1/(g-1)\leq m/n$, yielding a~contradiction.

\subsection{Thin fronts.} \label{ssect:tor4} $\,$ \vskip5pt

\nin Let us now describe the cells which will, in a~sense specified later, form a~core of the cubical complexes $TF_{m,n,g}$.

\begin{df} 
We call a torus front $T$ {\bf thin} if $w(T)\leq\sqrt 2$; we call it 
{\bf very thin} if $w(T)<\sqrt 2$.  
\end{df}

It is easy to construct a~very thin torus front for fixed $m$ and $n$: simply start from the point $(0,0)$ eastwards and then approximate the line connecting $(0,0)$ with $(m,n)$ as closely as possible, staying on one side of this line (touching the line is permitted). It is clear that if the approximation is not the closest possible, then $w(T)\geq\sqrt 2$. 

It is also important to notice that, since we assumed $m(g-1)< n$, each horizontal leg in a~thin torus front has length at least $g$, this can be proved with the same argument as the one we used to show that flipping all extreme corners preserves that property as well.

Let $I$ denote the line segment connecting $(0,0)$ with $(m,n)$. Counting from the origin in the northeastern direction, the first point on $I$ with integer coordinates, other than the origin itself, is $(m/\gcd(m,n),n/\gcd(m,n))$, which can be achieved in $(m+n)/\gcd(m,n)$ steps, (it is customary to denote the greatest common divisor with simple round brackets, but since $(m,n)$ is used to denote the point on the plane with these coordinates, we prefer here a~somewhat longer notation). 

The very thin torus front, which we just constructed, will have $(0,0)$ as a~northwestern extreme corner, and after that the northwestern extreme corners will repeat with the period of $(m+n)/\gcd(m,n)$. From this we see that in general, to determine a~very thin torus front we just need to specify, where on the part of the torus front connecting two consequent northwestern extreme corners the origin will lie. Therefore, we conclude that there are precisely $(m+n)/\gcd(m,n)$ very thin torus fronts, and that the width of any very thin $(m,n)$-torus front depends only on $m$~and~$n$.

We now extend the definition of width. For any flip $F$, we define its width $w(F)$ to be the maximum of widths of all the torus fronts which are members of~$F$. The concepts of millstones and extreme corners can be extended to flips as well, see Figure~\ref{fig:flip1}.

We call a~flip $F$ thin if $w(F)\leq\sqrt2$. The thin flips play
an~important role in our context, as the next theorem indicates.

\begin{thm} \label{thm:tfcoll}
For any natural numbers $m,n$, and $g$, there exists a~sequence of collapses leading from the cubical complex $TF_{m,n,g}$ to its subcomplex consisting of all thin flips.
\end{thm}

\pr To start with, we notice that if $\sigma$ and $\tau$ are cells of 
$TF_{m,n,g}$, and $\sigma$ is contained in $\tau$, then $w(\sigma)\leq
w(\tau)$. Therefore, it is enough to describe how to collapse away all
the cells of a~certain width ($>\sqrt 2$), under the condition that
the cells of larger width are already gone.

Let us fix this width and denote it by $w$. To every cell $\sigma$, we
can associate cells $\sigma^\uparrow$ and $\sigma_\downarrow$ as
follows: $\sigma^\uparrow$ is obtained from $\sigma$ by adding
elementary flips on all extreme corners which do not yet have them,
whereas $\sigma^\downarrow$ is obtained from $\sigma$ by removing the
elementary flips on all extreme corners, where in each elementary flip
we keep the extreme corner (since $w>\sqrt 2$ both corners of
an~elementary flip cannot be extreme, hence this is
well-defined). Note that it is possible that $\sigma=\sigma^\uparrow$,
or $\sigma=\sigma_\downarrow$, but not both.

By our discussion above, if all horizontal legs have length at least $g$ in vertices of $\sigma$, then the same is true for all cells in the interval $[\sigma_\downarrow,\sigma^\uparrow]$. We also see that, by our construction, all cells in the interval $[\sigma_\downarrow,\sigma^\uparrow]$ have the same width. Moreover, each such interval is isomorphic to a~Boolean algebra $B_t$, where $t$ is equal to the number of extreme corners. Clearly, the intervals $[\sigma_\downarrow,\sigma^\uparrow]$ are disjoint, and the removal of such an~interval constitutes a~collapse, as long as $\sigma^\uparrow$ is the only maximal cell containing $\sigma_\downarrow$ at the given time. 

This observation allows us to do as follows: restrict the partial inclusion order to the cells $\sigma^\uparrow$, choose a~linear extension of this partial order, and then  remove intervals $[\sigma_\downarrow,\sigma^\uparrow]$ following the linear extension in the decreasing order. The statement follows now from the fact, that if $\tau$ contains $\sigma$, then $\tau^\uparrow$ contains $\sigma^\uparrow$, since at the moment when the interval $[\sigma_\downarrow,\sigma^\uparrow]$ is being removed, this will translate to: if $\tau$ contains $\sigma_\downarrow$, then $\tau^\uparrow=\sigma^\uparrow$. 
\qed

\vskip5pt

Let us now finish our study of torus fronts by analyzing the
subcomplex $\thin_{m,n,g}$ of $TF_{m,n,g}$ consisting of all thin
flips. 

\begin{prop} \label{pr:thin}
For any natural numbers $m,n$, and $g$, such that $n>m(g-1)$, the
complex $\thin_{m,n,g}$ consists of $(m+n)/\gcd(m,n)$ cubes, each one
having dimension $\gcd(m,n)$, which are connected together by their
diagonally opposite end vertices, like a~garland, to form one cycle.
\end{prop}
\pr Let $v$ be a~thin torus front. Flipping a~corner which is not extreme
will give a~flip of width strictly larger than $\sqrt 2$, hence this
is not allowed. If, on the other hand, an~extreme northwestern corner
is flipped, then, the width of this flip is equal to $\sqrt 2$, and
flipping another extreme northwestern corner will not change the
width, whereas flipping any other vertex will increase it. The same
holds for flipping extreme southeastern corners.  It follows that
there are two ways to get a~maximal cell adjacent to $v$: either by
flipping all extreme northwestern corners, or by flipping all extreme
southeastern corners.

On the other hand, if $F$ is a~maximal cell, we have $w(F)=\sqrt 2$.
If we now choose a~torus front $T$ belonging to this cell, we will
also have $w(T)=\sqrt 2$, except for two cases: if in each elementary
flip we either always choose northwestern path, or if we always choose
southeastern path.  We see that each maximal cell contains exactly two
torus front of width strictly smaller than $\sqrt 2$, and that these
are opposite corners of the cube. Since $\thin_{m,n,g}$ is connected,
the conclusion follows.
\qed

\vskip5pt

We remark, that since each cube can further individually be collapsed
onto a~path connecting two opposite vertices, Theorem~\ref{thm:tfcoll}
and Proposition~\ref{pr:thin} entail that $TF_{m,n,g}$ can be
collapsed onto a~cycle.

\subsection{The implications for the cohomology groups of $\thom(C_m,K_n)$.} 
\label{ssect:tor5} $\,$ \vskip5pt

\nin  We now have enough information to complete the
computation of the entire second tableau, except for the $(m-1)$th
column, the latter containing the cohomology groups of
$\thom(C_m,K_n)$. Since the cohomology groups of $\thomp(C_m,K_n)$ are
non-trivial only in one dimension, we can then use this to derive
almost complete information about the cohomology groups of
$\thom(C_m,K_n)$ with $\zz$-coefficients. In fact, with a~little extra
effort we can describe $H^*(\thom(C_m,K_n);\zz)$ completely already
now.  For the sake of brevity we postpone the complete determination
to Section~\ref{s:zsect} where the general case of integer
coefficients will be dealt with.

By Propositions \ref{pr:red1} and \ref{pr:tm2}, to compute the values
$E_2^{p,q}$, for $p\neq m-1$, it is enough to calculate the
$\zz$-cohomology groups of the complexes $\Phi_{t,n,3}$, for the
required range of the parameter~$t$. These, on the other hand, are
isomorphic to $H^*(\thin_{t,n-t,3};\zz)$, this follows from
Proposition~\ref{pr:tfint}, and Theorem~\ref{thm:tfcoll}. Finally, by
Proposition~\ref{pr:thin}, we see that
$H^0(\thin_{t,n-t,3};\zz)=H^1(\thin_{t,n-t,3};\zz)=\zz$, and the other
cohomology groups are trivial. Tracing back all the indices we obtain
the following answer.

\begin{prop} \label{pr:e2z2}
The non-zero entries $E_2^{p,q}$, for $p\neq m-1$, are the following:
$E_2^{m-3,n-2}=\zz$, $E_2^{m-t-2,t(n-2)}=E_2^{m-t-1,t(n-2)}=\zz$, for
$2\leq t\leq\lfloor(m-1)/3\rfloor$, and
$E_2^{0,0}=\zz$. If~additionally $3\big{|}m$, then
$E_2^{2m/3-1,m(n-2)/3}={\mathbb Z}^3_2$.
\end{prop}

The non-zero entry $E_2^{2m/3-1,m(n-2)/3}={\mathbb Z}^3_2$ stems from
the fact that the complexes $\Phi_{m,mg,g}$ are disjoint unions of
$g$ points, and are not homotopy equivalent to circles.

An alert reader may also inquire why $E_2^{m-2,n-2}=0$. If the entries
in $(m-1)$th column were all $0$, we would of course have
$E_2^{m-2,n-2}=\zz$.  However, we know that
$H^{m+n-4}(\thomp(C_m,K_n);\zz)=0$, unless $(m,n)=(7,4)$. Since all
the entries to the northwest from $E_2^{m-2,n-2}$, which could
potentially eliminate this entry at a~later stage of the spectral
sequence computation, are $0$, we have no choice but to conclude that
the map
\[
d_1:E_1^{m-2,n-2}\big{/}\im(d_1:E_1^{m-3,n-2}\ra E_1^{m-2,n-2})\lra E_1^{m-1,n-2}
\]
is injective, and so $E_2^{m-2,n-2}=0$. The case $(m,n)=(7,4)$ can be
computed directly.

Since the reduced cohomology of $\thomp(C_m,K_n)$ is concentrated in
one dimension, Proposition \ref{pr:e2z2} allows us to compute almost
all cohomology groups of $\thom(C_m,K_n)$. The complete computation is
done in Section~\ref{s:zsect}, see Theorem~\ref{thm:main} for full
answer.

\section{Euler characteristic formula}

\nin While the cohomology groups, and even the Betti numbers of the 
$\thom$ complexes are usually very hard to compute, the Euler
characteristic may turn out to be a~more accessible invariant.

To start with, let $T$ and $G$ be arbitrary graphs. A~simple counting
in the filtration which we imposed on the simplicial complex
$\thomp(T,G)$ yields the following formula:
\begin{equation}\label{eq:eul1}
\wti\chi(\thomp(T,G))=\sum_{\emptyset\neq S\subseteq V(T)}
(-1)^{|S|+1}\wti\chi(\thom(T[S],G)),
\end{equation}
compare also with Proposition~\ref{pr:rho}.

\begin{thm} \label{thm:chi}
Let $T$ and $G$ be arbitrary graphs, then we have
\begin{equation} \label{eq:chi}
\wti\chi(\thom(T,G))=\sum_{\emptyset\neq S\subseteq V(T)}
(-1)^{|S|+1}\wti\chi(\thomp(T[S],G)).
\end{equation}
\end{thm}
\pr Take the formula \eqref{eq:eul1}, and apply M\"obius inversion 
on the Boolean algebra without the minimal element of all subsets of $V(T)$.
\qed

\begin{thm} \label{thm:chikn}
For an arbitrary graph $T$, we have the following formula:
\begin{equation} \label{eq:chikn}
\wti\chi(\thom(T,K_n))=\sum_{\emptyset\neq S\subseteq V(T)}
(-1)^{n+|S|}\wti\chi(\ind(T[S]))^n.
\end{equation}
\end{thm}

\pr By Corollary \ref{crl:hompkn} we know that 
$\thomp(T,K_n)=\ind(T)^{*n}$.  On the other hand, for any two
simplicial complexes $X$ and $Y$, we have
$\wti\chi(X*Y)=-\wti\chi(X)\wti\chi(Y)$. It follows that
$\wti\chi(\thomp(T,K_n))=(-1)^{n-1}\wti\chi(\ind(T))^n$. Substituting
this into identity \eqref{eq:chi} yields formula \eqref{eq:chikn}.
\qed

\vskip5pt

This gives a new proof of the following result.

\begin{crl} \label{crl:chi1} {\rm(\cite[Proposition 4.6]{BK03b},
\cite[Proposition 3.3.4]{IAS}).}

\nin For arbitrary positive integers $m$ and $n$, such that $n\geq m$, we
have
\begin{equation} \label{eq:chi1}
 \wti\chi(\thom(K_m,K_n))=\sum_{k=1}^{m-1} (-1)^{n+k+1}\binom{m}{k+1}k^n.
\end{equation}
\end{crl}

\pr Since $\ind(K_t)$ is just a~set of $t$ points, we have 
$\wti\chi(K_t)=t-1$.  Substituting this into equation \eqref{eq:chikn}
and noticing that the summands depend on the cardinality of $S$ only,
we obtain~\eqref{eq:chi1}.
\qed

\vskip5pt

If we try to use formula \eqref{eq:chikn} to compute the Euler
characteristic of $\thom(C_m,K_n)$, we will need nontrivial reductions
to get a~nice answer. Using Proposition~\ref{pr:e2z2} instead yields
a~very simple formula directly. This is reflecting very well the
reduction which often occurs when passing between tableaux in the
spectral sequence computation, since using the filtration on the
simplicial complex $\thomp(T,G)$ as we did is the same as using the
first tableau of our spectral sequence.

\begin{crl} \label{crl:chi2}
 For arbitrary positive integers $m$ and $n$, such that $m\geq 5$,
 $n\geq 4$, we have
\begin{equation} \label{eq:chi2}
 \wti\chi(\thom(C_m,K_n))=
\begin{cases}
(-1)^{nk-m},&\text{ if } m=3k\pm 1,\\
(-1)^{nk-m}(2^n-3),&\text{ if } m=3k.
\end{cases}
\end{equation}
\end{crl}

\pr The entries $E_2^{m-t-2,t(n-2)}$, and $E_2^{m-t-1,t(n-2)}$ cancel 
out in pairs, except for $E_2^{m-3,n-2}$, which we can discount for,
since $E_2^{m-2,n-2}$ already canceled out with $E_2^{m-1,n-2}$. Thus
the only nontrivial contribution comes from the Betti numbers of the
total complex $\thomp(C_m,K_n)$, and, in case $m=3k$ from the
entry $E_2^{2m/3-1,m(n-2)/3}$.
\qed

\section{Cohomology with integer coefficients} \label{s:zsect}

\subsection{Fixing orientations on $\thom$ and $\thomp$ complexes.} \label{ssect:z1} $\,$ 
\vskip5pt

\nin In this section we work only with integer coefficients, 
so we shall suppress $\dz$ from the notation. To be able to work with
integer coefficients we need to choose orientations on cells, so that
the differential maps can be determined. We shall use the notation and
we shall cite (without proofs) results from~\cite{BK03c}.

Let $T$ and $G$ be two graphs, and let us fix an~order of the vertices
of $T$ and of $G$. When convenient, we may identify vertices of~$T$,
resp.\ of $G$, with integers $1,\dots,|V(T)|$, resp.\
$1,\dots,|V(G)|$, according to the chosen orders. First we deal with
the simpler, simplicial case of $\thomp(T,G)$.  The vertices of
$\thomp(T,G)$ are indexed with pairs $(x,y)$, where $x\in V(T)$, $y\in
V(G)$, such that if $x$ is looped, then so is~$y$. We order these
pairs lexicographically: $(x_1,y_1)\prec(x_2,y_2)$ if either
$x_1<x_2$, or $x_1=x_2$ and $y_1<y_2$. Then, we orient each simplex of
$\thomp(T,G)$ according to this order on the vertices. We call this
orientation {\it standard}, and call the oriented simplex~$\eta_+$. We
identify this simplex with the corresponding chain, and denote the
dual cochain with~$\eta^*_+$.

Next, we deal with $C^*(\thom(T,G))$. On each cell $\eta\in\thom(T,G)$ we fix an~orientation (also called standard) as follows: orient each simplex $\eta(i)$ according to the chosen order on the vertices of $G$, then, order these simplices in the direct product according to the chosen order on the vertices of~$T$. We call this oriented cell $\eta$; a~choice of orders on the vertex sets of $T$ and $G$ is implicit.

Note, that permuting the vertices of the simplex $\eta(i)$ by some $\sigma\in\cs_{|\eta(i)|}$ changes the orientation of the cell $\eta$ by $\sgn\sigma$, and that swapping the simplices with vertex sets $\eta(i)$ and $\eta(i+1)$ in the direct product changes the orientation by $(-1)^{\dim\eta(i) \cdot\dim\eta(i+1)}$. Furthermore, if $\ti\eta\in \thom^{(i+1)}(T,G)$ is obtained from $\eta\in\thom^{(i)}(T,G)$ by adding a~vertex~$v$ to the list $\eta(t)$, then $[\eta:\ti\eta]$ is $(-1)^{k+d-1}$, where $k$ is the position of $v$ in $\ti\eta(t)$, and $d$ is the dimension of the product of the simplices with the vertex sets $\eta(1),\dots,\eta(t-1)$, i.e., $d=1-t+\sum_{j=1}^{t-1}|\eta(j)|$. 

Let $\wti C^*$ be the subcomplex of $C^*(\thomp(T,G))$ generated by
all $\eta^*_+$, for $\eta:V(T)\ra 2^{V(G)}$, such that
$\supp\eta=V(T)$.  Set
\begin{equation} \label{eq:x*}
X^*(T,G):=\wti C^*[|V(T)|-1].
\end{equation}
Note that both $C^i(\thom(T,G))$ and $X^i(T,G)$ are free $\dz$-modules
with the bases $\{\eta^*\}_\eta$ and $\{\eta^*_+\}_\eta$ indexed by
$\eta:V(T)\ra 2^{V(G)}\sm\{\emptyset\}$, such that $\sum_{j=1}^{|V(T)|}|\eta(j)|=|V(T)|+i$. Furthermore, when $\eta$ is a~cell of $\thom(T,G)$, we set 
$$c(\eta):=\sum_{\begin{subarray}{c}i\text{ is even}\\
1\leq i\leq|V(T)|\end{subarray}}
|\eta(i)|.$$
For any $\eta:V(T)\ra 2^{V(G)}\sm\{\emptyset\}$, set $\rho(\eta_+):= (-1)^{c(\eta)}\eta$. The induced map $\rho^*:X^i(T,G)\ra C^i(\thom(T,G))$ is of course an~isomorphism of abelian groups for any~$i$. It turns out that more is true.

\begin{prop} \label{pr:rho}
{\rm (\cite[Proposition 3.3]{BK03c}).}

\nin  For any two graphs $T$, and $G$, the map $\rho^*:X^*(T,G)\ra C^*(\thom(T,G))$ 
is an~isomorphism of the cochain complexes.
\end{prop}

\subsection{Signed versions of formulas for generators $[\sigma^S_V]$.} \label{ssect:z2} $\,$ \vskip5pt

\nin Let us now return to considering the cohomology classes $[\sigma^S_V]$ defined in Subsection~\ref{ssect:arc2}, though this time we are working with integer coefficients. Clearly, these still index the generators of $E_1^{p,t(n-2)}$. The new feature appearing when working over integers is that the sign of the generator may change when vertex representatives move along edges. 

The pattern of the sign change is the same as in the case $|V|=1$
described in \cite{BK03c}. Indeed, let $S$ be a~proper subset of
$[m]$, so that $C_m[S]$ is a~forest. Let $(S,A)$ be an~arc picture,
and let $V$ be the set of vertex representatives. Choose an~arc $a\in
A$, and denote its vertex representative by $v\in a$. Choose $w\in a$, 
such that $(v,w)\in E(C_m)$, and set $W:=(V\cup\{w\})\sm\{v\}$. We would
like to understand the sign change which occurs when $v$ is replaced
by~$w$ as the vertex representative of~$a$. Let $\coprod_{A} K_2$
denote the disjoint union of copies of $K_2$ indexed by $A$. We define
a~graph homomorphism $\coprod_{A} K_2\hookrightarrow C_m$ as follows:
for the copy of $K_2$ indexed by $a$ we choose one of the two graph
homomorphisms $K_2\ra(v,w)$, for the other copies of $K_2$ we take any
graph homomorphism mapping one of the vertices of $K_2$ to the vertex
representative of the corresponding arc, and mapping the other vertex
of $K_2$ to any of the neighboring vertices in~$S$.

This homomorphism induces a~$\zz$-equivariant map
$H^*(\thom(\coprod_{A} K_2,K_n))\ra H^*(\thom(C_m[S],K_n))$, where the
$\zz$-action on the left hand side is induced by the antipodal action
on the copy of $K_2$ indexed by~$a$.  On the other hand, we have
$\thom(\coprod_{A}K_2,K_n)\cong\prod_{A}\thom(K_2,K_n)\cong\prod_{A}S^{n-2}$ and the $\zz$-action on the last space is the antipodal action on the factor indexed by~$a$.  Recall, that the antipodal action on $S^n$ changes the sign
of the $n$-dimensional cohomology generator if and only if $n$ is
even, and the same is true if we are acting on one of the factors in
a~direct product. Using already introduced terminology, we can take
$\rho^*(\sigma_V^S)$ as a~representative for the generator of
$H^{t(n-2)}(\thom(C_m[S],K_n))$, which we are considering. If
$(v,w)\neq(1,m)$, then there is no further sign change and we get
\[[\rho^*(\sigma_V^S)]=(-1)^{n+1}[\rho^*(\sigma_W^S)].\]

If on the other hand $(v,w)=(1,m)$, then the $\zz$-action does not interchange two neighboring simplices, but the first one and the last one instead. One of these simplices has dimension 0 and the other one has dimension $n-2$. Using the rules for the sign change described in the subsection \ref{ssect:z1}, we see that we get additional sign factor $(-1)^{n(|A|-1)}$, since we have to swap simplices of dimension $n-2$ in total $|A|-1$ times, and each swap yields the sign $(-1)^{(n-2)^2}=(-1)^n$; we ignore the swaps involving 0-dimensional simplices since they do not influence the sign. So in this case we get
\[[\rho^*(\sigma_V^S)]=(-1)^{n|A|+1}[\rho^*(\sigma_W^S)].\]

Combining with Proposition~\ref{pr:rho} we get
\begin{equation} \label{eq:shift1}
  \left[\sigma^S_V\right]=-\left[\sigma^S_W\right],
\end{equation} 
if $(v,w)\neq(1,m)$, since the swapped vertices have different parity. In case $(v,w)=(1,m)$, the swapped vertices may or may not have different parity, so we have an~additional sign factor $(-1)^{nb}$, where $b$ is the number of vertices in $S$ between $w$ and~$v$. So in this case we get 
\begin{equation} \label{eq:shift2}
  \left[\sigma^S_V\right]=(-1)^{n|A|+nb+n+1}\left[\sigma^S_W\right].
\end{equation} 

Next, we look at the analog of the formulas \eqref{eq:d1gen} and \cite[(4.7)]{BK03c}. Combining the construction of generators $\sigma_V^S$ with the choice of orientations on $\thomp(T,K_n)$, we get the following formula.

\begin{equation} \label{eq:d1genz}
d_1(\left[\sigma_{V}^S\right])=
\sum_{w\notin S}(-1)^{\sgn(w)}\left[\sigma_{V}^{S\cup\{w\}}\right],	 
\end{equation}
where the sum is taken over all $w$, such that adding $w$ to $S$ does not unite two arcs from $A$, and the sign $(-1)^{\sgn(w)}$ is given by the formula \begin{equation} \label{eq:sgnw}
\sgn(w)=|S\cap[w-1]|+n\cdot|V\cap[w-1]|.
\end{equation}
The same way as for $\zz$-coefficients, the equation \eqref{eq:d1genz} is also valid when $|S|=m-1$.

\subsection{Completing the calculation of the second tableau.} \label{ssect:z3} $\,$ \vskip5pt

\nin We shall now do the same computation as we did in Subsection \ref{ssect:tor5}, only this time with integer coefficients.

The first obstruction is that Proposition \ref{pr:tm2} may not be valid anymore, since we have to take the signs into account. However, though the cochain complexes $(A^*_t,d_1)$ and $(C_*(\Phi_{t,m,3};\dz),d)$ may not be isomorphic, the only difference is that some of the incidence numbers differ by a~sign.

The argument, using torus fronts, with which we computed $H_*(\Phi_{t,m,3};\zz)$, consisted of presenting a~sequence of collapses leading from $TF_{t,m-t,3}$ to $\thin_{t,m-t,3}$. Due to isomorphism of cochain complexes over $\zz$, stated in Proposition \ref{pr:tm2}, these collapses could have been performed directly on $(A_t^*,d_1)$. 

It is now crucial to realize, that the same collapses can still be performed in the cochain complex $(A_t^*,d_1)$, even though we are working over integers. This follows from the version of Discrete Morse theory for chain complexes, see \cite{dmt}, since our matching is acyclic, and the weights on the matching relations, which here are the incidence numbers, are $\pm 1$, hence invertible over~$\dz$. 

Let $m-t>2t$, i.e., $m>3t$, and consider the complexes $\thin_{t,m-t,3}$. On each maximal cube, choose any of the shortest paths connecting the opposite vertices by which these cubes hang together. Such a~path has $\gcd(t,m-t)=\gcd(t,m)$ edges. It is obvious that the collapsing process can be continued until the whole complex $\thin_{t,m-2,3}$ is collapsed onto this path~$P$. Let us extend this collapsing to $(A_t^*,d_1)$ as well, and let us denote the remaining chain complex by $(Q_t,d_1)$.

The chain complex $(Q_t,d_1)$ has only two non-zero entries, so let us write it as $0\ra Q_t^1\ra Q_t^0\ra 0$. If we worked over $\zz$, this chain complex would be simply computing the homology of a~circle with coefficients in~$\zz$. Over integers we have signs and there are two possibilities: either $H_1(Q_t)=H_0(Q_t)=\dz$, or $H_1(Q_t)=0$, and $H_0(Q_t)=\zz$. To distinguish between these two cases we shall now calculate the differential. 

Both groups $Q_t^1$ and $Q_t^0$ have $m$ generators. The boundary of each generator of $Q_t^1$ is a~sum of exactly two generators of $Q_t^0$ with coefficients $\pm 1$. Let $\alpha$ be the number of those generators of $Q_t^1$ where these two coefficients are equal. It is easy to see that $H_1(Q_t)=H_0(Q_t)=\dz$, if $\alpha$ is even, and $H_1(Q_t)=0$, $H_0(Q_t)=\zz$, if $\alpha$ is odd, so all that remains is to calculate the parity of~$\alpha$.

\begin{figure}[hbt]
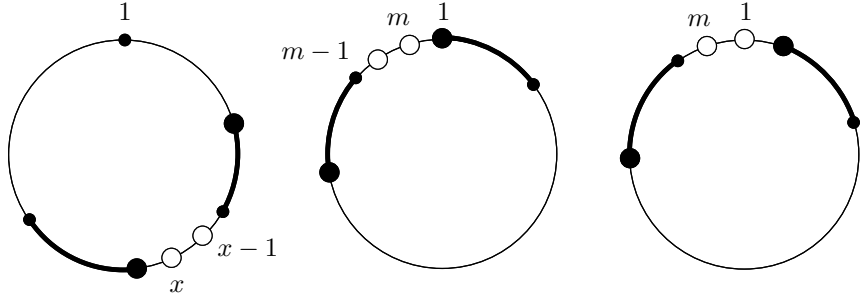

\begin{center}
  \begin{picture}(0,0)%
    \includegraphics{3cases.pstex}%
  \end{picture}%
  \input{3cases.pstex_t}%
  
\end{center}
\caption{3 cases of the sign calculation.}
\label{fig:3cases}
\end{figure}

One may think of generators of $Q_t^1$ as being indexed by $t$-arc
pictures $(S,A)$, with all but one gaps consisting of one element, and
one 2-element gap, denoted $[x-1,x]_m$. Taking the boundary
corresponds to narrowing this gap in two different ways: by extending
one of the bordering arcs clockwise or counterclockwise. We break our
calculation into 3 cases. See Figure~\ref{fig:3cases}.

\nin {\it Case 1.} $2\leq x\leq m-1$. By formula \eqref{eq:d1genz}, filling in $x-1$ or $x$ will have the same sign. Furthermore, when filling in $x$, we shall also need to move the vertex representative from $x+1$ to $x$. These vertices have different parity in the order on $V(C_m[S\cup\{a\}])$, hence by \eqref{eq:shift1} we get an additional factor $-1$ for this generator. We conclude that all these generators of $Q_t^1$ do not contribute to~$\alpha$. 

\nin {\it Case 2.} $x=m$. Again, by formula \eqref{eq:d1genz}, filling in $x-1$ or $x$ will have the same sign. However, in this case we have to be careful about the sign caused by moving the vertex representative from $1$ to $m$. By formula \eqref{eq:shift2}, the additional sign is $(-1)^{mn+n+1}$, since $b=m-t-2$, and $|A|=t$. So the contribution to $\alpha$ is $1$ if $mn+n$ is odd, and $0$ otherwise.

\nin {\it Case 3.} $x=1$. In this case shifting the vertex representative gives the sign $-1$, by \eqref{eq:shift1}. The non-trivial contribution comes instead from the fact that filling $x$ and $x-1$ may give different signs. The difference in signs is $(-1)^{nt+m+t+1}$, since $m-t-1+t(n-2)$ is total number of vertices (if counted with multiplicities) in positions $2$ through $m-1$. Therefore, in this case the contribution to $\alpha$ is $1$ if $nt+m+t+1$ is even, and $0$ otherwise.

\vspace{5pt}

By construction of the path $P$, the group $Q_t^1$ contains exactly one generator with $x=m$, and exactly one generator with $x=1$. Summarizing we see that $\alpha$ is odd if and only if $mn+n+nt+m+t+1$ is odd.  Note that $mn+n+nt+m+t+1=(m+t+1)(n+1)$. Thus we obtain the extension of Proposition \ref{pr:e2z2} to the case of integer coefficients. We present it in the Table~\ref{tab:zans}.

\begin{table}[hbt]
\[
\begin{array}
{|l|c|c||} 
\hline
&& \\
n,m,t \,\,\, & \,\,\,\,\,\,E_2^{m-t-2,t(n-2)}\,\,\,\,\,\, &
\,\,\, \,\,\,E_2^{m-t-1,t(n-2)}\,\,\, \,\,\,\\ 
&& \\ \hline
t=1  && \\
m(n+1) \text{ is odd, }  & 0 & 0\\
m(n+1) \text{ is even, } & \dz &0\\ 
\hline
n \text{ is odd, }&&\\
\text{or } m+t \text{ is odd, } & \dz &\dz\\ 
2\leq t\leq\lfloor (m-1)/3\rfloor  && \\
\hline
n \text{ is even, }&&\\
\text{and } m+t \text{ is even, } & 0 &\zz\\ 
2\leq t\leq\lfloor (m-1)/3\rfloor  && \\
\hline
n \text{ is any, }&&\\
\text{and } 3\big{|}m, & 0 &\dz^3\\ 
t=m/3  && \\
\hline

\hline
\end{array}
\]
\caption{$\,$}
\label{tab:zans}
\end{table}

\subsection{Summary: the full description of the groups $\wti H^*(\thom(C_m,K_n);\dz)$.} 
\label{ssect:z6} $\,$ \vskip5pt

\nin For an abelian group $\Gamma$, we let $\Gamma(d)$ denote the copy of
$\Gamma$ in a~graded $\dz$-algebra, placed in dimension~$d$.

\begin{thm}\label{thm:main}
For any integers $m,n$, such that $m\geq 5$, $n\geq 4$, we have 
\begin{equation}\label{eq:main}
\wti H^*(\thom(C_m,K_n);\dz)=\left(\bigoplus_{t=1}^{\lfloor(m-2)/3\rfloor}
A_{t,m,n}\right)
\oplus B_{m,n},
\end{equation}
where
\begin{equation}\label{eq:main1}
A_{t,m,n}=\begin{cases}
\dz(tn-3t)\oplus\dz(tn-3t+1),& \text{ if } n \text{ is odd or } m+t \text{ is odd,}\\
\zz(tn-3t+1),& \text{ if } n \text{ is even and } m+t \text{ is even,}
\end{cases}
\end{equation}
and
\begin{equation}\label{eq:main2}
B_{m,n}=\begin{cases}
\dz^{2^n-3}(nk-m),& \text{ if } m=3k,\\
\dz(nk-m+2),& \text{ if } m=3k+1,\\
\dz(nk-m),& \text{ if } m=3k-1.
\end{cases}
\end{equation}
\end{thm}

For example $\wti H^*(\thom(C_6,K_4);\dz)=A_{1,6,4}\oplus
B_{6,4}=\dz(1)\oplus\dz(2)\oplus\dz^{13}(2)=\dz(1)\oplus\dz^{14}(2)$,
and $\wti H^*(\thom(C_8,K_6);\dz)=A_{1,8,6}\oplus A_{2,8,6}\oplus B_{8,6}=
\dz(3)\oplus\dz(4)\oplus\zz(7)\oplus\dz(10)$.

 \vskip5pt

\nin {\bf Proof of Theorem~\ref{thm:main}.} 
Recall that all cohomology groups of $\thomp(C_m,K_n)$ are zero,
except for one, see~Corollary \ref{crl:cplus}. Since our spectral
sequence converges to $H^*(\thomp(C_m,K_n);\dz)$, we know that almost
all entries in the second tableau should cancel out.

Since $n\geq 4$, the calculation summarized in Table~\ref{tab:zans}
implies that there will be no non-zero differential between entries in
columns $0,\dots,m-2$.  If $n\geq 5$, then the differentials with
entries in the $(m-1)$st column as a~target will never have the same
target entry. If $n=4$, the targets for two such differentials may be
the same, but then the sources will simply sum up. This follows from
the pattern of the entries and from the algebraic fact that if we have
a~group homomorphism $f:A\ra G$, such that $f$ is an~injection, $A$ is
either $\dz$ or $\zz$, and $G/\im f$ is isomorphic to $\dz$, then $G$
is isomorphic to $A\oplus\dz$.

Since $\wti H^i(\thomp(C_m,K_n))=0$, for $i<nk-1$, we obtain the ``A
part'' of~\eqref{eq:main}, with the exception of a~single entry in the
case $m=3k+1$, which is dealt with below. So we are done with the
calculation of almost all cohomology groups of $\thom(C_m,K_n)$ at
this point. To get the remaining ``B part'', we need to see what
happens on the diagonals $x+y=const\geq nk-2$. We shall consider 3
different cases.

Assume first that $m=3k-1$. We see from Table~\ref{tab:zans} that the
top (in terms of the sum of coordinates) non-zero element is
$E_2^{m-k,(n-2)(k-1)}$. This element is on the diagonal
$m-k+(n-2)(k-1)=m+nk-3k-n+2=nk-n+1=(nk-1)-(n-2)$, hence it is $n-2$
diagonals away from the diagonal $nk-1$, and we conclude that
$B_{m,n}=\dz(nk-m)$.

Assume now that $m=3k$. The top non-zero element in columns
$0,\dots,m-2$ is $E_2^{m-k-1,k(n-2)}$, which lies precisely on the
diagonal $x+y=nk-1$. By the analysis in Subsection \ref{ssect:arc2},
the generators of $E_2^{m-k-1,k(n-2)}$ are given by those
$[\sigma_{A_\bu}^S]$, where $S$ is obtained from $[m]$ by deleting
every third element, which can be done in three different ways, and
$A$ is the collection of all $m/3$ arcs of length~$2$.  By
construction, we see that in this case $d(\sigma_{A_\bu}^S)=0$, where
the differential is taken in $C^*(\thomp(C_m,K_n))$.  Additionally,
there are no non-zero elements northwest from $E_2^{m-k-1,k(n-2)}$. It
follows that the spectral sequence differential
$d_k:E_2^{m-k-1,k(n-2)}\ra E_2^{m-1,nk-m+1}$ is a~zero-map, and
therefore $E_2^{m-k-1,k(n-2)}=\dz^3$. Choosing a~field $\cf$ of
arbitrary characteristic we see that the Betti number of the entry
$E_2^{m-1,nk-m}$ is equal to $2^n-3$ if $n\geq 5$, and $2^n-2$, if
$n=4$. Since this does not depend on the choice of $\cf$, we conclude
that $H^{nk-m}(\thom(C_m,K_n))=E_2^{m-1,nk-m}$ is equal to the direct
sum if the corresponding number of copies of $\dz$, and hence
$B_{m,n}=\dz^{2^n-3}(nk-m)$ in this case.

Finally, assume that $m=3k+1$. In this case we have to deal with the
differentials from $E_2^{m-k-2,k(n-2)}=E_2^{m-k-1,k(n-2)}=\dz$, and,
if $n=4$, also from $E_2^{m-k,2k-2}=\zz$. Let $X_+$ denote
$\thomp(C_m,K_n)$, let $X\subset X_+$ be the subcomplex consisting of
all cells $\eta$, such that $|\supp\eta|<m$, and consider the
cohomology long exact sequence of the pair $(X_+,X)$. Since the
reduced cohomology groups $\wti H^*(X_+;\dz)$ are trivial in all
dimensions except for $nk-1$, we find the following exact sequence inside
the considered one:
\begin{multline}\label{eq:lesx}
0\longleftarrow H^{nk}(X_+,X)\longleftarrow H^{nk-1}(X)\stackrel f\longleftarrow \\
\stackrel f\longleftarrow H^{nk-1}(X_+)\longleftarrow H^{nk-1}(X_+,X)\longleftarrow H^{nk-2}(X)\longleftarrow 0.
\end{multline}
The crucial observation now if that the map $f$ is an~inclusion map to
a~direct summand. To see this, set
$Y:=\thomp(C_m[V(C_m)\sm\{1\}],K_n)\subset
X$. By~\cite[Proposition~5.2]{K2}, and by the fact that
$\thomp(C_m,K_n)\simeq\ind(C_m)^{*n}$, we see that the inclusion map
$i:Y\hookrightarrow X_+$ induces isomorphism of the cohomology groups
$i^*:H^*(X_+)\ra H^*(Y)$. The isomorphism map $i^{nk-1}$ factors as
$H^{nk-1}(X_+)\stackrel f\ra H^{nk-1}(X)\ra H^{nk-1}(Y)$, since $f$
itself is induced by the inclusion map $X\hookrightarrow X^+$. It
follows from this factorization that $f$ is an~inclusion map to
a~direct summand, so the sequence~\eqref{eq:lesx} splits in this
sense.

Reinterpreting this for our original spectral sequence we see that
this means that the differential $d_{k+1}:E_2^{m-k-2,k(n-2)}\ra
E_2^{m-1,nk-m+1}$ is a~zero-map, and that, if $n=4$, then the
differential $d_{k-1}:E_2^{m-k,2k-2}\ra E_2^{m-1,k}$ is
an~isomorphism. This proves the ``B part'' as well as the rest of the
``A part'' in this case.
\qed

\vskip5pt

We remark that in general the homotopy type of $\thom(C_m,K_n)$ is
unknown, even conjecturally.

\section{Grinding complexes of maps between cycles}

\nin In this section we take another look at the complexes of graph
homomorphisms between cycles $\thom(C_m,C_n)$. For having case-free
statements, we assume that $m\geq 5$, $n\geq 3$, $n\neq 4$, the other
cases are very simple. It was shown in~\cite{CK1} that the homotopy
type of each connected component of $\thom(C_m,C_n)$ is either that of
a~point, or that of~$S^1$. The complexes $\thom(C_m,C_n)$ are cubical,
and the encoding of cells, suggested in \cite[Section 5]{CK1} can be
easily translated into the language of torus fronts as follows. We
rephrase Definition~\ref{df:torfr} introducing an extra parameter.

\begin{df} \label{df:torfr2}
Let $b$ be a~positive integer, and consider the action of the group
$\dz\times\dz$ on the plane, where the standard generators of
$\dz\times\dz$ act by vector translations, with vectors $(-b,b)$ and
$(m,n)$. A {\bf $b$-band $(m,n)$-torus front} (or sometimes simply
band torus front) is an~orbit of this $(\dz\times\dz)$-action on the set
of all $(m,n)$-grid paths.
\end{df} 

As an~analogous generalization, the complexes $TF_{m,n,g}^b$ are
defined by replacing $(m,n)$-torus fronts in
Definition~\ref{df:tfcomp} by $b$-band $(m,n)$-torus fronts.

It is clear that, when translated into our language, the description
from \cite[Section 5]{CK1} would simply say that the cells of
$\thom(C_m,C_n)$ are encoded by band torus fronts. More precisely,
consider the subcomplex $X_w$ of $\thom(C_m,C_n)$ induced by the
vertices which have the winding number $w$ (each vertex is a~graph
homomorphism between cycles, which gives rise to continuous map of
$S^1$ to itself, hence the notion of winding number is well defined).
Without loss of generality we can assume that $w$ is
nonnegative. Then, $m-nw$ is even, and we set $r=(m-nw)/2$. If both
$m$ and $n$ are even, the complex $X_w$ has two connected components,
otherwise it is connected. Tracing the constructions we obtain the
following reformulation.

\begin{prop}
The cubical complex $X_w$ is isomorphic to the band torus front
complex $TF_{r,m-r,1}^n$.
\end{prop}

Most of the discussion from Section \ref{sect5}, including such
notions as extreme corners, thin flips, and millstones, generalizes to
the band torus front context without any modification whatsoever. For
the sake of brevity, we do not repeat it here, and simply state the
main result.

\begin{thm} \label{thm:tfcoll2}
For any natural numbers $m,n$, and $g$, there exists a~sequence of
collapses leading from the cubical complex $TF_{m,n,g}^b$ to its
subcomplex consisting of all thin $b$-band flips.
\end{thm}

We now obtain the band analog of Proposition~\ref{pr:thin}.

\begin{crl}\label{crl:cmcn}
When $r>0$, the cubical complex $X_w$ collapses to the subcomplex
consisting of $mn/\gcd(m,r)$ cubes of dimension $\gcd(m,r)$. These
cubes are linked by the opposite vertices to form $2$ circles of both
$m$ and $n$ are even, and one circle otherwise.
\end{crl}
\pr The proof of Proposition~\ref{pr:thin} holds. The only difference is that when enumerating the cubes we have an~extra factor $n$ to compensate for working on an~$n$-band.
\qed

\vskip5pt

The cubes in the thin part of $\thom(C_m,C_n)$ are linked in vertices,
which encode those graph homomorphisms from $C_m$ to $C_n$ which
distribute their turning points as evenly as possible.

Since collapsing $\thom(C_m,C_n)$ onto the subcomplex of thin graph
homomorphisms produces a~$(\cd_m\times\cd_n)$-equivariant strong
deformation retraction, Corollary~\ref{crl:cmcn} yields in particular
a~full description of the $(\cd_m\times\cd_n)$-action on the
cohomology groups of $\thom(C_m,C_n)$, where $\cd_t$ denotes the
$t$-th dihedral group.

\vskip10pt

\nin {\bf Acknowledgments.} We would like to thank the Swiss National 
Science Foundation and ETH-Z\"urich for the financial support of this
research.

\end{document}